\documentclass[11pt,a4paper]{article} 
\usepackage[latin1]{inputenc}
\usepackage[T1]{fontenc}
\usepackage{lmodern}
\usepackage[final]{changes}
\usepackage{varioref}
\usepackage{indentfirst}
\usepackage{ulem}
\usepackage[english]{babel}
\usepackage{amsmath}
\usepackage{amsthm}
\usepackage{amsfonts}
\usepackage[vmargin=2.5cm,hmargin=3.3cm]{geometry}
\usepackage{amssymb}
\usepackage{pstricks-add}
\usepackage{graphicx} 
\usepackage{color}
\usepackage{wallpaper}
\usepackage[english]{minitoc}
\usepackage{fancyhdr} 
\color{black}
\usepackage{fancyhdr}
\usepackage{url}
\usepackage{hyperref} 
\hypersetup{
colorlinks=true, 
breaklinks=true, 
urlcolor= blue, 
linkcolor= red, 
citecolor=green, 
pdftitle={}, 
pdfauthor={Lahoucine Elaissaoui}, 
pdfsubject={Simulation} 
}

\def\N {\mathbb{N}}

\def\d {\mathrm{d}}

\newcommand{\equa}[1]{\begin{equation} #1 \end{equation}}

\newtheorem{thee}{\sc Theorem}
\newtheorem{propo}{Proposition}[section]
\labelformat{propo}{\sc Proposition~#1}
\newtheorem{lemm}{\sc Lemma}
\labelformat{lemm}{\sc Lemma~#1}
\newtheorem{defin}{Definition}[section]

\newtheorem{coro}{\sc Corollary}
\labelformat{thee}{\sc Theorem~#1}
\labelformat{coro}{\sc Corollary~#1}
\title{  \textsc{Evaluation of Log-tangent Integrals by series involving $\zeta(2n+1)$  } \vspace*{0.3 cm} }
\author{\textsc{BY} \\ \textsc{Lahoucine Elaissaoui}\footnote{\texttt{E-mail:lahoumaths@gmail.com}} \qquad \ \quad \textsc{And} \quad \textsc{Zine El Abidine Guennoun}\footnote{E-mail: guennoun@fsr.ac.ma} \\ \\ \\ \sc{Mohammed V University in Rabat} \\ \sc{Faculty of Sciences} \\ \sc{Department of Mathematics}  \\  \\ \sc{Morocco}  }

\date{}
\begin{document}
\maketitle
\begin{center}
\textit{This is not a final version! But, the Original and peer-reviewed version is published by Taylor \& Francis Group in Integral Transforms and Special Functions. Is available on 07/04/2017 \texttt{http://www.tandfonline.com/10.1080/10652469.2017.1312366}.}
\end{center}
\vspace*{0.3 cm}

\begin{abstract}
In this note, we show that the values of integrals of the log-tangent function with respect to any square-integrable function on $\left[0 , \frac{\pi}{2} \right]$ may be determined (or approximated) by  an infinite (or finite) sum involving the Riemann Zeta-function at odd positive integers. 
\end{abstract}

\vspace*{0.5 cm}

{\bf Key Words:} Riemann Zeta function, Ap\'{e}ry's constant, Catalan's constant, Summation formula, Log-tangent integrals, Finite and Infinite series.

{\bf 2010 Mathematics Subject Classification (s):} 11M06, 26D15, 11L03.

\vspace*{1 cm}

\section{Introduction and preliminaries}

 The Riemann Zeta-function, denoted $\zeta$, is defined by
 $$\zeta(s) := \sum_{k = 1}^{\infty} \frac{1}{k^s}. $$
 \replaced{The}{ the} sum \added{$\zeta(s)$} is absolutely convergent for any complex number in the half-plane $\Re s > 1$ and \deleted{the Riemann Zeta-function} is analytic on this half-plane. In particular, if $s$ is a positive integer greater than $1$ it is well-known  that
  $$ \zeta(s=2n) = \frac{(-1)^{n-1} 2^{2n}B_{2n}}{2(2n)!} \pi^{2n} \qquad \replaced{(n \in \mathbb{N})}{ \forall n \geq 1},$$ 
  where $B_n$ denotes the $n$-th Bernoulli number. \added{Here and in the following, let $\mathbb{R}$ and $\mathbb{N}$ be the sets of real numbers and positive integers, respectively, and let $\mathbb{N}_0 = \mathbb{N}\cup\left\{ 0\right\}$}. For the odd numbers, i.e. $s=2n+1$, no closed forms have been proven yet. However, it has been conjectured by Kohnen~\cite{Koh} that the ratios of the quantities $\frac{\zeta(2n+1)}{\pi^{2n+1}}$ are transcendental for every integer \replaced{$n \in \mathbb{N}$}{ $n \geq 1$}. 
 
  Ap\'{e}ry's constant is defined as the number 
 $$\zeta(3) = 1.202056903159594285399738161511449990764986292 .... $$
 It was named for the French mathematician Roger Ap\'{e}ry who proved in 1978~\cite{Ape} that it is irrational; Ap\'{e}ry's theorem. However, it is still not known whether Ap\'{e}ry's constant is transcendental. Recently, T. Rivoal~\cite{Riv} and W. Zudilin~\cite{Zud} \replaced{have}{ has} shown, respectively, that infinitely many of the numbers $\zeta(2n+1)$ must be irrational, and that at least one of the eight numbers \replaced{$\zeta(2n+1)$ $(n=2,\cdots, 9)$}{$\zeta(5)$, $\zeta(7)$,..., $\zeta(17)$, $\zeta(19)$} must be irrational.
  
 Of course, several series and integrals involving the numbers $\zeta(2n+1)$, with $n \geq 1$, have been shown; see \added{,} for example \added{,} \cite{Sri0,Con,Choi2,choi,choi33,choi4,Sri1,Sri2,Sri}. As we notice, almost all these results were obtained by evaluation of log-sine integrals and its related functions. In fact, the log-sine integrals were firstly introduced by Euler in 1769 and he showed in 1772~\cite{Eul} that 
 $$ \int_0^{\frac{\pi}{2}} x \log\left( \sin x\right) \d x = \frac{7}{16} \zeta(3)  - \frac{\pi^2}{8} \log 2;$$

and by exploiting this \replaced{integral}{ integrals}, Euler gave the famous series representation of Ap\'{e}ry's constant \replaced{as follows:}{ , given below}

$$ \zeta(3) = \frac{\pi^2}{7} \left( 1 - 2 \sum_{n=1}^{\infty} \frac{\zeta(2n)}{2^{2n}(2n+1)(n+1)}\right).$$
Notice that there are many other classical series representation of Ap\'{e}ry's constant. The higher moments were studied differently by several authors (see\added{,} for example \added{,} \cite{Koy}). Furthermore, many families of log-sine and log-cosine integrals were evaluated explicitly by Choi and Srivastava; see \cite{Choi1} and \cite{choi3}.

In this note, we shall study integrals involving the log-tangent function for a certain class of functions $f$ defined on the interval $\left[0, \frac{\pi}{2} \right]$; namely
\equa{L(f) := \int_0^{\frac{\pi}{2}} f(x) \log \left( \tan x \right) \d x. \label{L}}
Moreover, we show that these integrals (for some class of functions) may be approximated by a series of terms involving the numbers $\zeta(2n+1)$, with \replaced{$n \in \mathbb{N}$}{$n \in \mathbb{N}^*$}. 

In fact, the integral in \eqref{L} exists whenever the function belongs to $L^2\left( \left[0, \frac{\pi}{2}\right]\right)$ and by Cauchy-Schwarz inequality we have
$$ \frac{2}{\pi}\left|L(f)\right| \leq \frac{\pi}{2} \left\Vert f \right\Vert_2,$$ 
where 
$$\Vert f \Vert_2 = \sqrt{\frac{2}{\pi}\int_0^{\frac{\pi}{2}}\left[ f(x)\right]^2 \d x};$$
notice that the constant on the right-hand side of inequality above follows from the quantity \cite[eq. 5, p. 533]{Jef.Zwi}
\equa{\frac{2}{\pi}\int_0^{\frac{\pi}{2}}\left( \log\left( \tan x\right)\right)^2 \d x = \frac{\pi^2}{4}. \label{LTN}}
Also, if $f$ is bounded on $\left[ 0 , \frac{\pi}{2} \right]$ then we have 
$$\left|L(f) \right| \leq 2G \Vert f \Vert_{\infty},  $$
where $G$ is Catalan's constant, 

$$G = - \int_0^{\frac{\pi}{4}} \log\left( \tan x\right) \d x = \sum_{k \geq 0} \frac{(-1)^k}{(2k+1)^2} ,$$
and $\Vert \cdot \Vert_{\infty}$ is the supremum norm. Notice that, the table \cite{Jef.Zwi} contains explicit evaluations of the integral in \eqref{L} for some  trigonometric functions $f$ (see p. 533-591). However, the case when $f$ is polynomial, which is strongly connected with the numbers $\zeta(2n+1)$, is not treated; for example, if $f(x) = x$ and $f(x) = x^2$ we have the following

   \begin{propo}
   Ap\'{e}ry's constant can be represented \deleted{as follows} by
    $$ \int_0^{\frac{\pi}{2}} x \log \left( \tan x\right) \d x = \frac{7}{8} \zeta(3) $$
    and 
    $$\int_0^{\frac{\pi}{2}} x^2 \log \left( \tan x\right) \d x = \frac{7}{16}\pi \zeta(3). $$
    \label{Th}\end{propo}
    \textbf{Proof.} Bradley showed in \cite[Th.1]{Bra} that for $x \in \left[0, \frac{\pi}{2} \right]$
    \equa{ \int_0^x \log \left( \tan u\right) \d u = - \sum_{n=0}^{\infty} \frac{\sin\left(2(2n+1)x \right)}{(2n+1)^2}; \label{Fs}}
    notice that, the series in the right-hand side is absolutely convergent for all \replaced{$x \in \mathbb{R}$}{ $x \in \left[0, \frac{\pi}{2} \right]$}. Then
    $$\int_0^{\frac{\pi}{2}} \left[\int_0^x \log \left( \tan u\right) \d u\right] \d x = - \sum_{n \geq 0}\frac{1}{(2n+1)^3}.$$
   We apply integration by parts for the integral in the left-hand side, we obtain
   $$\int_0^{\frac{\pi}{2}} \left[\int_0^x \log \left( \tan u\right) \d u \right] \d x= - \int_0^{\frac{\pi}{2}} x \log \left( \tan x \right) \d x. $$
   Since, for all $s>1$, 
    $$\sum_{n \geq 0} \frac{1}{(2n+1)^s} = \left( 1 - \frac{1}{2^{s}}\right) \zeta(s) $$
   then 
       $$ \int_0^{\frac{\pi}{2}} x \log \left( \tan x \right) \d x = \frac{7}{8} \zeta(3). $$
       
     Similar reasoning, using the Fourier expansion \eqref{Fs}, yields
     $$ \int_0^{\frac{\pi}{2}} x^2 \log \left( \tan x\right) \d x = - \int_0^{\frac{\pi}{2}}2x \left[\int_0^x \log\left( \tan u \right)\d u\right] \d x = \frac{7}{16}\pi \zeta(3).$$ 
     \begin{flushright}
      $\square$
      \end{flushright} 

 The next section contains the evaluation of integrals $L(P)$ for any given polynomial $P$. Thereby, we will deduce that any square-integrable function may be approximated or determined by a sum involving $\zeta(2n+1)$. We conclude, in the last section, with a brief discussion of corresponding results with some additional remarks. 
 
  \section{Evaluation of the integral $L$ for polynomial functions}
 
  Let us start with the following
   
  \begin{lemm}
  For any positive integer $k$ we have
  $$\int_0^{\frac{\pi}{2}} \cos(2kx) \log\left( \tan x\right)\d x = \begin{cases}\begin{array}{cll} 0 & \text{if} &  k \ \text{is even} \\ -\frac{\pi}{2k} & \text{if}& k \ \text{is odd} \end{array} \end{cases} .$$
 \label{Lm1} \end{lemm}
  \textbf{Proof.} If $k=2n$, then
  \begin{align*}
  \int_0^{\frac{\pi}{2}} \cos(4nx) \log\left( \tan x\right)\d x &= \int_0^{\frac{\pi}{2}}\cos\left( 4n \left( \frac{\pi}{2}-x\right)\right) \log \left( \tan \left( \frac{\pi}{2} - x \right) \right) \d x \\ &= - \int_0^{\frac{\pi}{2}}\cos(4nx) \log \left( \tan x \right) \d x .
  \end{align*}
  Therefore,
  $$ \int_0^{\frac{\pi}{2}} \cos(4nx) \log\left( \tan x\right)\d x = 0 .$$
  Now, for $k=2n+1$, using integration by parts and the Fourier expansion \eqref{Fs}, we have
  \begin{align*}\int_0^{\frac{\pi}{2}} \cos(2kx) \log\left( \tan x\right)\d x &= 2(2n+1) \int_0^{\frac{\pi}{2}} \sin\left(2(2n+1)x\right) \left[ \int_0^x \log\left( \tan u \right) \d u \right] \d x \\ &= -(2n+1) \sum_{m=0}^{\infty}\frac{1}{(2m+1)^2}\int_0^{\frac{\pi}{2}}\left[ \cos\left( 4(m-n) x\right) - \cos \left(4(m+n)x \right) \right]\d x \\ &= - \frac{\pi}{2} \frac{1}{2n+1} ;
  \end{align*}
which completes the proof.
\begin{flushright}
$\square$
\end{flushright}

Now, we shall generalize results of \ref{Th} for higher moments, namely
 $$\int_0^{\frac{\pi}{2}} x^n \log\left( \tan x\right) \d x, \qquad (n \in \N^*). $$
 
  For this fact, let us briefly recall some properties of Euler polynomials, denoted $E_n$. It is well-known that Euler polynomials are defined on the unit interval $[0,1]$ and they are Appell sequences. Moreover, \replaced{Euler polynomials $E_n(x)$ are defined by the following generating function}{ the generating function for Euler polynomials is defined by}
 
  \equa{\frac{2e^{tx}}{e^t+1} = \sum_{n \geq 0} E_n(x) \frac{t^n}{n!}, \qquad |t|<\pi .\label{GeF}}
  Of course, the polynomials $E_n$ have several interesting properties, the most important to us are given below:
  \begin{itemize}
  \item Symmetry; for all $x$ in $[0,1]$ 
  \equa{E_n(1-x) = (-1)^nE_n(x), \ \replaced{(n \in \mathbb{N}_0)}{ \qquad n\geq 0}. \label{Sym}}
  \item Inversion; every monomial $x^n$ may be expressed in terms of $E_n(x)$, namely
  \equa{x^n = E_n(x) + \frac12 \sum_{k=0}^{n-1}\binom{n}{k} E_k(x).\label{Inv}}
  for every integer $n > 0$ and for all $x  \in [0,1]$.
  \item Translation; for all $x \in [0,1]$ and a given real $y$ we have
  \equa{ E_n(x+y) = \sum_{k=0}^n\binom{n}{k} E_k(x)y^{n-k}.\label{Tra}}
  \item The Fourier series form for the Euler polynomials for a positive integer $n \geq 1$: 
  \equa{ E_{2n-1}(x) = 4(-1)^n(2n-1)!\pi^{-2n} C_n(x)\label{FrS}} 
  and
  $$ E_{2n} = 4(-1)^n(2n)! \pi^{-2n-1}S_n(x),$$ 
  for all $x$ in the unit interval; the functions $C_n$ and $S_n$ are defined by
  $$C_n(x) := \sum_{k \geq 0}\frac{\cos\left((2k+1) \pi x\right)}{(2k+1)^{2n}} $$
  and
  $$S_n(x) := \sum_{k \geq 0}\frac{\sin\left((2k+1)\pi x\right)}{(2k+1)^{2n+1}}  .$$ 
  \end{itemize}
Now, we can easily prove the following
\begin{thee}
Let $n$ be a positive integer, then 
$$\int_0^{\frac{\pi}{2}}E_{2n}\left( \frac{2}{\pi}x\right) \log\left( \tan x \right) \d x = 0 $$
and
\equa{\int_0^{\frac{\pi}{2}}E_{2n-1}\left( \frac{2}{\pi}x\right) \log\left( \tan x \right) \d x = \frac{(-1)^{n-1} (2n-1)!}{\pi^{2n-1}}(2-2^{-2n}) \zeta(2n+1). \label{Eulodd}}
\label{Th1}\end{thee}
\textbf{Proof.} The first integral is immediate from the property \eqref{Sym}. 
For the second integral we use formula \eqref{FrS}, hence
\begin{align*} \int_0^{\frac{\pi}{2}}E_{2n-1}\left( \frac{2}{\pi}x\right) \log\left( \tan x \right) \d x &= 4(-1)^n(2n-1)!\pi^{-2n} \int_0^{\frac{\pi}{2}}C_n\left( \frac{2}{\pi}x\right) \log\left( \tan x\right) \d x \\ &=\frac{4(-1)^n(2n-1)!}{\pi^{2n}} \sum_{k=0}^{\infty} \frac{1}{(2k+1)^{2n}} \int_0^{\frac{\pi}{2}}\cos\left(2(2k+1)x\right) \log \left( \tan x\right) \d x.
\end{align*}
Then \ref{Lm1} completes the proof of \ref{Th1}, namely

\begin{align*}  \int_0^{\frac{\pi}{2}}E_{2n-1}\left( \frac{2}{\pi}x\right) \log\left( \tan x \right) \d x &= \frac{2(-1)^{n-1}(2n-1)!}{\pi^{2n-1}} \sum_{k=0}^{\infty} \frac{1}{(2k+1)^{2n+1}} \\ &= \frac{2(-1)^{n-1}(2n-1)!}{\pi^{2n-1}} \left( 1- \frac{1}{2^{2n+1}}\right)\zeta(2n+1).
\end{align*}

\begin{flushright}
  $\square$
  \end{flushright}  
  
Notice that, we can deduce from formula \eqref{Eulodd}, for any integer $n \geq 1$, that
$$ \int_0^{\frac{\pi}{4}} E_{2n-1}\left(\frac{2}{\pi} x \right) \log\left( \tan x\right) \d x = \frac{(-1)^{n-1} (2n-1)!}{\pi^{2n-1}}\left(1-\frac{1}{2^{2n+1}}\right) \zeta(2n+1).$$
Moreover, the inversion formula \eqref{Inv} implies the following
 
  \begin{coro}
Let $n$ be a positive integer, then we have
  \begin{align*}
  \int_0^{\frac{\pi}{2}} x^n \log \left( \tan x \right) \d x &= (-1)^{\lfloor \frac{n-1}{2}\rfloor}\frac{n!}{2^{n-1}}\left(1-\frac{1}{2^{n+2}}\right)\zeta(n+2) \delta(n) \\ &+ \frac{n!}{2^{n}} \sum_{k=1}^{\lfloor{\frac{n}{2}}\rfloor} \frac{(-1)^{k-1}\pi^{n-2k+1}}{(n-2k+1)!} \left(1-\frac{1}{2^{2k+1}}  \right) \zeta(2k+1), 
  \end{align*}
  where $\lfloor \cdot \rfloor$ is the floor function and 
  $$\delta(n) := \begin{cases}\begin{array}{cll} 1 & \text{if} & n \ \text{is odd} \\ 0 & \text{if} & n \ \text{is even.}\end{array} \end{cases}$$
  
   \label{Ct}\end{coro}
 
 One can also use the translation property \eqref{Tra} to show that
 \begin{coro} For any positive integer $n$ and a given real $y$,
 $$\int_0^{\frac{\pi}{2}}E_{2n}\left( \frac{2}{\pi}x + y \right) \log\left( \tan x\right) \d x = (2n)! \sum_{k=1}^n \frac{(-1)^{k-1}}{\pi^{2k-1}}\left( 1 - \frac{1}{2^{2k+1}}\right)\zeta(2k+1) \frac{y^{2n-2k+1}}{(2n-2k+1)!}, $$
 and
 $$\int_0^{\frac{\pi}{2}}E_{2n-1}\left( \frac{2}{\pi}x + y \right) \log\left( \tan x\right) \d x = (2n-1)! \sum_{k=1}^n \frac{(-1)^{k-1}}{\pi^{2k-1}}\left( 1 - \frac{1}{2^{2k+1}}\right)\zeta(2k+1) \frac{y^{2n-2k}}{(2n-2k)!}. $$
 \end{coro}
 
 Actually, we can extract more identites involving the integrals $L(E_n)$ by using other formulas for Euler polynomials. Moreover, similar reasoning to that used to obtain the results above is applicable for other similar polynomials such as the Bernoulli polynomials.  Of course, \ref{Ct} permits us to evaluate more integrals $L(P)$ for any given polynomial $P$; an explicit evaluation of $L(P)$ for a polynomial $P$ is given in the following
 
 \begin{thee} Let $P$ be polynomial of degree \replaced{$m \in \mathbb{N}$. Then}{ $m \geq 1$, then we have}
 
 $$\int_0^{\frac{\pi}{2}}P(x) \log \left( \tan x\right)\d x = \sum_{k=1}^{\lfloor \frac{m+1}{2} \rfloor} \frac{(-1)^{k-1}}{2^{2k-1}}\left[P^{(2k-1)}\left( \frac{\pi}{2}\right) + P^{(2k-1)}(0) \right] \left( 1 - \frac{1}{2^{2k+1}}\right)\zeta(2k+1). $$
 where $P^{(p)}(\alpha)$ denotes the $p$-th derivative of $P$ at the point $\alpha$.
 \label{Th2}
 \end{thee}
  
  \textbf{Proof.} It is easy to show, using integration by parts, that
  $$ \int_0^{\frac{\pi}{2}}P'(x) \sin\left( 2(2n+1)x\right) \d x = \frac{P'\left( \frac{\pi}{2} \right)+P'(0)}{2(2n+1)} - \frac{1}{2^2(2n+1)^2} \int_0^{\frac{\pi}{2}}P^{(3)}(x) \sin\left( 2(2n+1)x\right) \d x.$$
  Thus by induction we obtain 
  $$ \int_0^{\frac{\pi}{2}}P'(x) \sin\left( 2(2n+1)x\right) \d x = \sum_{k=1}^{\lfloor \frac{m+1}{2}\rfloor} \frac{(-1)^{k-1}}{2^{2k-1}}\left[ P^{(2k-1)}\left( \frac{\pi}{2}\right) + P^{(2k-1)}(0) \right] \frac{1}{(2n+1)^{2k-1}}.$$
  Since, by formula \eqref{Fs},
  \begin{align*}
  L(P) &= - \int_0^{\frac{\pi}{2}} P'(x) \left[\int_0^{x} \log\left( \tan u \right)\d u \right] \d x \\ &= \sum_{n = 0}^{\infty} \frac{1}{(2n+1)^2} \int_0^{\frac{\pi}{2}}P'(x) \sin\left( 2(2n+1)x\right) \d x;  
  \end{align*}
  then 
  $$L(P) = \sum_{k=1}^{\lfloor \frac{m+1}{2}\rfloor} \frac{(-1)^{k-1}}{2^{2k-1}}\left[ P^{(2k-1)}\left( \frac{\pi}{2}\right) + P^{(2k-1)}(0) \right]\sum_{n = 0}^{\infty} \frac{1}{(2n+1)^{2k+1}} .$$
  Thereby, the fact that
  $$ \sum_{n = 0}^{\infty} \frac{1}{(2n+1)^{2k+1}} = \left( 1 - \frac{1}{2^{2k+1}}\right) \zeta(2k+1)$$
  completes the proof.
  \begin{flushright}
  $\square$
  \end{flushright}
  
  There are many applications of \ref{Th2}, the most interesting one is that it allows to expand the log-tangent function in a series of orthogonal polynomials corresponding to Hilbert space $L^2\left(0 , \frac{\pi}{2} \right)$ with the inner product
  $$\langle f , g \rangle := \frac{2}{\pi}\int_0^{\frac{\pi}{2}}f(x)g(x) \d x , \qquad \forall f,g \in L^2\left(0 , \frac{\pi}{2} \right).$$  
  We define on $\left[0, \frac{\pi}{2}\right]$ the shifted Legendre polynomials  $\left(\tilde{P}_n\right)_{n \geq 0}$ using Rodrigues formula
  $$\tilde{P}_n(x) := \frac{(-1)^n}{n!}\left( \frac{2}{\pi} \right)^n \frac{\d^n}{\d x^n}\left[ x^n\left( \frac{\pi}{2}- x\right)^n\right].  $$
Indeed, the shifting function $x \mapsto \frac{4}{\pi}x -1$ bijectively maps the interval $\left[0 , \frac{\pi}{2} \right]$ to the interval $[-1,1]$ in which the ordinary Legendre polynomials $P_n$ are defined (see for example \cite{Sam} for more details about Legendre polynomials). Therefore, the shifted Legendre polynomials form an orthogonal basis in $L^2\left(0 , \frac{\pi}{2} \right)$ and we have
  $$\langle \tilde{P}_n , \tilde{P}_m\rangle = \begin{cases} \begin{array}{cll} \frac{1}{2n+1} & \text{if} & n=m \added{,}\\ 0 & \text{if}& n \neq m \added{.}\end{array} \end{cases} \deleted{.} $$
Notice that the polynomials $\tilde{P}_n$ satisfy two following properties
\begin{itemize}
\item For any given integer $n \geq 0$ and all $x \in \left[ 0 , \frac{\pi}{2}\right]$ we have
\equa{\tilde{P}_n\left( \frac{\pi}{2} - x \right) = (-1)^n \tilde{P}_n(x). \label{SymL} }
\item An explicit representation is given by
$$\tilde{P}_n(x) = (-1)^n \sum_{k=0}^n \binom{n}{k}\binom{n+k}{k} \left( - \frac{2}{\pi}\right)^k x^k. $$
\end{itemize}
The last property implies that the $m$-th derivative of the shifted Legendre polynomials is 
$$\tilde{P}_n^{(m)}(x) = (-1)^n m! \sum_{k=m}^n \binom{n}{k}\binom{n+k}{k} \binom{k}{m} \left( - \frac{2}{\pi}\right)^k x^{k-m};$$
thereby, for any integer $ 0 \leq m \leq n$ we have
\equa{ \tilde{P}_n^{(m)}(0) = (-1)^{n+m} \left(\frac{2}{\pi}\right)^m \frac{(n+m)!}{m!(n-m)!}. \label{SLP0}}
Consequently, the coefficients of the log-tangent function in the orthogonal basis $\left\{\tilde{P}_n\right\}_{n \geq 0}$ of Hilbert space $L^2\left( 0 , \frac{\pi}{2} \right)$ are given in the following

\begin{coro} \replaced{For}{ for} any positive integer $n$\added{, }we have
$$L\left( \tilde{P}_{2n} \right) = 0, $$
and
$$ L\left( \tilde{P}_{2n-1} \right) = 2\sum_{k=1}^{n}\frac{(-1)^{k-1}}{\pi^{2k-1}} \frac{(2(n+k-1))!}{(2k-1)!(2(n-k))!}\left( 1 - \frac{1}{2^{2k+1}}\right) \zeta(2k+1) . $$
\label{CoeffLT}\end{coro}
  \textbf{Proof.} The first integral follows directly, using integration by substitution, utilizing property \eqref{SymL}. 
 The second integral follows by applying \ref{Th2} and the fact that 
 $$ \tilde{P}_{2n-1}^{(2k-1)}\left( \frac{\pi}{2} \right) =  \tilde{P}_{2n-1}^{(2k-1)}(0) = \left(\frac{2}{\pi}\right)^{2k-1} \frac{(2n+2k-2)!}{(2k-1)!(2n-2k)!}. $$
 Notice that the first and the second equalities follow respectively from \eqref{SymL} and \eqref{SLP0}.
  \begin{flushright}
  $\square$
  \end{flushright}

Furthermore, it is well-known that every element of Hilbert space $L^2$ can be written in a unique way as a sum of multiples of these base elements. Namely, for our case,
$$ \forall f \in L^2\left( 0 , \frac{\pi}{2}\right) \qquad f = \sum_{n=0}^{\infty} c_n\left( f \right) \tilde{P}_n,  $$
where
$$c_n\left(f\right) = \frac{\langle f, \tilde{P}_n\rangle }{\Vert \tilde{P}_n \Vert_2^2} = (2n+1)\frac{2}{\pi} \int_0^{\frac{\pi}{2}} f(x) \tilde{P}_n(x) \d x . $$
Moreover, we have Parseval's identity
$$\Vert f \Vert_2^2 = \sum_{n=0}^{\infty}\frac{\left| c_n\left( f \right)\right|^2}{2n+1} .$$

The equality \eqref{LTN} implies that log-tangent function belongs to Hilbert space $L^2\left( 0 , \frac{\pi}{2}\right)$; hence, we deduce the following results
\begin{coro} We have,
$$ \log\left( \tan x \right)\overset{L^2\left(0 , \frac{\pi}{2} \right)}{=} \frac{2}{\pi}\sum_{n=1}^{\infty}(4n-1)L\left( \tilde{P}_{2n-1}\right) \tilde{P}_{2n-1}(x),  $$
and
$$ \sum_{n = 1}^{\infty} (4n-1) \left[ L\left(\tilde{P}_{2n-1} \right)\right]^2 = \left( \frac{\pi}{2} \right)^4. $$
\label{LTexp}\end{coro} 

We should not forget to mention that log-tangent function is also defined by the series \cite[1.518, eq. 3, p. 53]{Jef.Zwi}
$$\log\left( \tan x \right)= \log x + 2 \sum_{k=1}^{\infty} \frac{2^{2k-1} - 1}{k} \zeta(2k) \left( \frac{x}{\pi}\right)^{2k}, \quad x\in \left( 0, \frac{\pi}{2} \right). $$
Therefore, one can extract more series representations involving the numbers $\zeta(2n+1)$ by combining the different results obtained in this paper with the series above. Moreover, one can use the log-tangent expansion showed in \ref{LTexp} to prove several identities involving the numbers $\zeta(2n+1)$ and other constants. For example, a series representaion of Catalan's constant 

\begin{align*}
G &= - \int_0^{\frac{\pi}{4}} \log\left( \tan x\right) \d x \\ &= \sum_{n=1}^{\infty}(4n-1)L\left( \tilde{P}_{2n-1}\right) \left( -\frac{2}{\pi} \int_0^{\frac{\pi}{4}}\tilde{P}_{2n-1}(x) \d x \right) \\ &= \sum_{n = 1}^{\infty} \frac{(-1)^{n-1}(4n-1)}{2^{2n}n} \binom{2n-2}{n-1}L\left( \tilde{P}_{2n-1} \right).
\end{align*}
Notice that the evaluation of the integral in the second line is due to Byerly \cite[p.172]{Bye}, namely (in particular)
$$ -\frac{2}{\pi} \int_0^{\frac{\pi}{4}}\tilde{P}_0(x)\tilde{P}_{2n-1}(x) \d x = (-1)^{n-1}\frac{(2n-2)!}{2^{2n}n!(n-1)!}=\frac{(-1)^{n-1}}{2^{2n}n} \binom{2n-2}{n-1}. $$ 

On the other hand, any function $f \in L^2\left( \left[ 0, \frac{\pi}{2} \right]\right)$ may be expanded in terms of shifted Legendre polynomials and it may be approximated by its partial sum
$$f_N = \sum_{n=0}^N c_n(f) \tilde{P}_n ,$$
for a sufficiently large integer $N$. Hence, the integral $L(f)$, for any square-integrable function on $\left[0, \frac{\pi}{2}\right]$, may be approximated by the partial sum $L(f_N)$ which depends on the numbers $\zeta(2n+1)$, with $n \geq 1$. Namely we have the following
  \begin{coro}
  \replaced{For}{ for} any square-integrable function $f$ on $\left( 0 , \frac{\pi}{2}\right)$, there exists a sequence of real numbers $\left\{c_{N,k}(f) \right\}$, where $N>0$ is a sufficiently large integer and $k=1,\cdots, N$,  such that
  $$L(f) = \lim_{N \to +\infty} \sum_{k=1}^N c_{N,k}(f) \frac{(-1)^{k-1}}{\pi^{2k-1}}\left(1 - \frac{1}{2^{2k+1}} \right) \zeta(2k+1).$$
  Futhermore,
  $$c_{N,k}(f) = 2\sum_{j=k}^{N} (4j-1) \langle f, \tilde{P}_{2j-1}\rangle \frac{(2j+2k-2)!}{(2k-1)!(2j-2k)!}.$$
  \label{Apprx}\end{coro}
  \textbf{Proof.} Let $f$ be a square-integrable function, then $f$ may be expanded as
  $$f(x) := \sum_{j=0}^{\infty} (2j+1)\langle f , \tilde{P}_j\rangle \tilde{P}_{j}(x) .$$
  Then
  \begin{align*}
  L(f) &= \frac{\pi}{2} \langle f , \log\tan\rangle \\ &= \sum_{j=1}^{\infty} (4j-1)\langle f , \tilde{P}_{2j-1}\rangle L\left( \tilde{P}_{2j-1}\right) \\ &= \lim_{N \to \infty} \sum_{j=1}^{N} (4j-1)\langle f , \tilde{P}_{2j-1}\rangle L\left( \tilde{P}_{2j-1}\right) \\ &=\lim_{N\to \infty} 2\sum_{j=1}^{N} \sum_{k=1}^{j}(4j-1)\langle f , \tilde{P}_{2j-1}\rangle \frac{(-1)^{k-1}}{\pi^{2k-1}} \frac{(2j+2k-2)!}{(2k-1)!(2j-2k)! } \left(1 - \frac{1}{2^{2k+1}} \right) \zeta(2k+1) \\ &= \lim_{N \to \infty} \sum_{k=1}^N \left(2\sum_{j=k}^N (4j-1)\langle f , \tilde{P}_{2j-1}\rangle \frac{(2j+2k-2)!}{(2k-1)!(2j-2k)! } \right)\frac{(-1)^{k-1}}{\pi^{2k-1}}\left(1 - \frac{1}{2^{2k+1}} \right) \zeta(2k+1).
  \end{align*}
  Therefore
  $$L(f) = \lim_{N \to +\infty} \sum_{k=1}^N c_{N,k}(f) \frac{(-1)^{k-1}}{\pi^{2k-1}}\left(1 - \frac{1}{2^{2k+1}} \right) \zeta(2k+1).  $$
  \begin{flushright}
  $\square$
  \end{flushright}
  Notice that the convergence of the sum
  $$ \sum_{k=1}^N c_{N,k}(f) \frac{(-1)^{k-1}}{\pi^{2k-1}}\left(1 - \frac{1}{2^{2k+1}} \right) \zeta(2k+1) $$
  is not always \replaced{uniform}{ uniforme}; however \ref{Apprx} provide us a good approximation. For example, the expansion of the function $f(x) = \sqrt{x}$ for $N=5$ is
$$f_5(x) = \frac{\sqrt{2\pi}}{3}\tilde{P}_{0}(x) + \frac{\sqrt{2\pi}}{15}\tilde{P}_{1}(x) - \frac{\sqrt{2\pi}}{105}\tilde{P}_{2}(x)+\frac{\sqrt{2\pi}}{315}\tilde{P}_{3}(x) - \frac{\sqrt{2\pi}}{693}\tilde{P}_{4}(x)+ \frac{\sqrt{2\pi}}{1287}\tilde{P}_{5}(x); $$
then, after \added{a} simplification, we find
\begin{align*} 
L(f_5) &= \frac{42}{13\sqrt{2\pi}}\zeta(3) - \frac{1581}{13 \pi^2 \sqrt{2\pi}}\zeta(5) + \frac{13335}{13 \pi^4 \sqrt{2\pi}}\zeta(7) \\ &\approx 0.688084888082269488....  
\end{align*}
However, using \texttt{Mathematica} we find that
$$ L(f) := \int_0^{\frac{\pi}{2}} \sqrt{x} \log\left( \tan x\right) \d x \approx 0.689247.$$

Actually, $L$ is a linear functional on $L^2\left(0 , \frac{\pi}{2} \right)$ then $L$ is surjective onto the scalar field $\mathbb{R}$. Consequently we obtain the following density result
\begin{coro}
For any real number $\alpha$ there exist a sequences of real numbers $\{c_{N,k}\}$; where $N \in \N^*$ and $k=1, \cdots , N$, such that
$$\alpha = \lim_{N \to +\infty} \sum_{k=1}^N c_{N,k} \frac{(-1)^{k-1}}{\pi^{2k-1}}\left(1 - \frac{1}{2^{2k+1}} \right) \zeta(2k+1). $$ 
\end{coro}
We should not forget to mention that a similar result above has been showed by Alkan \cite[Th.1]{Alk}.

It is well-known that any continuous function $f$ on $\left[0 , \frac{\pi}{2} \right]$ is integrable and square-integrable. Therefore, one can use a similar manipulation to evaluate the integral $L(f)$. Moreover, the Weierstrass approximation theorem and \ref{Th2} allow us to state that: for any continuous real-valued function $f$ on $\left[0 , \frac{\pi}{2} \right]$ there exists a sequence of real numbers $\{a_{n,k}\}_{k=1,\cdots,n}$,  with $n \geq 1$, such that
$$ L(f) = \lim_{n \to \infty} \sum_{k=1}^n (-1)^{k-1} a_{n,k} \left( 1 - \frac{1}{2^{2k+1}}\right) \zeta(2k+1). $$ 

\section{Concluding results and remarks}

  Exploiting the different results obtained in this note, one can evaluate further integrals involving the log-tangent function with respect to some non-polynomial functions, as the following example shows:
Let $z$ be a complex number such that $|z|<1$ and let $f(x,z) = e^{2zx}$, then by generating function formula \eqref{GeF} and \ref{Th1} we obtain 
  
\begin{align*}
\int_0^{\frac{\pi}{2}} e^{2zx} \log\left( \tan x \right) \d x &= \int_0^{\frac{\pi}{2}} e^{\left(\frac{2}{\pi}x\right)\pi z} \log\left( \tan x \right) \d x \\ &= \left( e^{\pi z} + 1\right) \sum_{n=0}^{\infty} \frac{\pi^n z^n}{n!} \int_0^{\frac{\pi}{2}}E_n\left( \frac{2}{\pi}x\right) \log\left( \tan x\right) \d x \\&= \left( e^{\pi z}+1 \right) \sum_{n=1}^{\infty} (-1)^{n-1} \left( 1 - \frac{1}{2^{2n+1}} \right)\zeta(2n+1) z^{2n-1}.
\end{align*}
Moreover, since
\begin{align*}
 \int_0^{\frac{\pi}{2}} e^{2zx} \log\left( \tan x \right) \d x &= \int_0^{\frac{\pi}{2}} e^{2z\left(\frac{\pi}{2} - x \right)} \log\left( \tan\left( \frac{\pi}{2} -  x\right) \right) \d x \\ &= - e^{\pi z} \int_0^{\frac{\pi}{2}} e^{-2zx} \log\left( \tan x \right) \d x, 
 \end{align*}
 then for all $|z|<1$,
 
 $$ \int_0^{\frac{\pi}{2}} \sinh\left( 2xz - \frac{\pi}{2} z \right) \log\left( \tan x \right) \d x = 2\cosh\left(\frac{\pi}{2} z \right) \sum_{n=1}^{\infty} (-1)^{n-1} \left( 1 - \frac{1}{2^{2n+1}} \right)\zeta(2n+1) z^{2n-1}.$$
Also, we have
\begin{align*}
\int_0^{\frac{\pi}{2}} \cos\left( 2xz  \right) \log\left( \tan x \right) \d x &= \frac12 \int_0^{\frac{\pi}{2}} \left[e^{2izx}+e^{-2ixz}\right] \log\left( \tan x \right) \d x \\ &= -\sin\left( \pi z \right) \sum_{n=1}^{\infty} \left( 1 - \frac{1}{2^{2n+1}} \right)\zeta(2n+1) z^{2n-1},
\end{align*}
on the other hand, using \ref{Lm1} and sine identity $\sin(a)\sin(b) = \frac{1}{2}\left[\cos\left( \frac{a-b}{2}\right) - \cos\left( \frac{a+b}{2}\right) \right]$, one can show that
$$\int_0^{\frac{\pi}{2}} \cos\left( 2xz  \right) \log\left( \tan x \right) \d x = \frac{\sin\left(  \pi z \right)}{4z} \left[ \psi\left( \frac{1+z}{2}\right) + \psi\left( \frac{1-z}{2}\right) - 2\psi\left( \frac12 \right) \right], $$
where $\psi$ is the digamma function and $\psi\left( \frac{1}{2} \right) = -\gamma - 2 \log 2$; here $\gamma \approx 0.57721...$ is the Euler-Mascheroni constant. We would like to mention that the partial sum $$S_N(z) := \sum_{n=1}^{N} \left( 1 - \frac{1}{2^{2n+1}} \right)\zeta(2n+1) z^{2n-1} $$ 
converges very quickly to the function $F(z):= -\frac{1}{4z} \left[ \psi\left( \frac{1+z}{2}\right) + \psi\left( \frac{1-z}{2}\right) - 2\psi\left( \frac12 \right) \right]$ for every complex $|z|<1$.

More generally, using \added{a} similar reasoning as in the proof of \ref{Th2} one can evaluate the integral $L$ for a particular class of functions. In fact, let $f$ be a smooth -real or complex valued- function defined on $\left[0 , \frac{\pi}{2}\right]$ such that     
 $$b_k := \sup_{0 \leq x \leq \frac{\pi}{2}}\left| f^{(k)}(x)\right| = o\left( 2^k\right),\qquad  \text{as} \ k \to \infty .$$ Then we have
  
  $$ \int_0^{\frac{\pi}{2}}f(x) \log\left( \tan x \right) \d x =\sum_{k = 1}^{\infty}\frac{(-1)^{k-1}}{2^{2k-1}}\left[f^{(2k-1)}(0) + f^{(2k-1)}\left(\frac{\pi}{2} \right) \right] \left(1-\frac{1}{2^{2k-1}}\right)\zeta(2k+1).  $$
 Or\added{,} alternatively, for a given $|z|\leq 1$,
 
$$ \int_0^{\frac{\pi}{2}}f(zx) \log\left( \tan x \right) \d x =\sum_{k = 1}^{\infty}\frac{(-1)^{k-1}}{2^{2k-1}}\left[f^{(2k-1)}(0) + f^{(2k-1)}\left(\frac{\pi}{2}z \right) \right] \left(1-\frac{1}{2^{2k-1}}\right)\zeta(2k+1) z^{2k-1}.$$
 
 It should be noted that the coefficients $L\left( \tilde{P}_{2n-1}\right)$ obtained in \ref{CoeffLT} converge very quickly to $0$ as $n \to +\infty$. Furthermore, the partial sum   
$$ \sum_{n = 1}^{N} \frac{(-1)^{n-1}(4n-1)}{2^{2n}n} \binom{2n-2}{n-1}L\left( \tilde{P}_{2n-1} \right),$$
converges very quickly to Catalan's constant as well; for instance, for $N=10$ 
$$\sum_{n = 1}^{10} \frac{(-1)^{n-1}(4n-1)}{2^{2n}n} \binom{2n-2}{n-1}L\left( \tilde{P}_{2n-1} \right)= 0.914611602803... \approx G - 1.3539\times10^{-3}.  $$

 Finally, we would like to note that several integrals given earlier  may play an important role, in number theory, in regard to the proof or disproof of the algebraicity of numbers in the form $\frac{\zeta(2n+1)}{\pi^{2n+1}}$. In addition, the constants $\zeta(2n+1)$  arise naturally in a number of physical and mathematical problems. For instance, in physics, they appear in correlation functions of antiferromagnetic xxx spin chain and when evaluating the $2n$-dimensional form of the Stefan-Boltzmann law.

\section*{Acknowledgements}
The authors are very grateful to Mr. James Arathoon for the endless English corrections. Furthermore, the authors would like to express their gratitude to the anonymous referees.

 \bibliographystyle{plain}  
  
\end{document}